\documentstyle[12pt,epsfig]{article}
\def\be{\begin{equation}}
\def\ee{\end{equation}}
\def\disp{\displaystyle}

\def\scri{\scriptsize}
\def\ve{\varepsilon}
\def\th{\theta}
\def\lam{\lambda}

\def\R{{\sf I\kern-.2em R}}
\def\C{\kern.1em{\raise.47ex\hbox{$\scriptscriptstyle |$}}\kern-.40em{\sf C}}
\def\Z{{\sf Z\kern-.32em Z}}

\newtheorem{theorem}{\noindent Theorem}
\newtheorem{lemma}{\noindent Lemma}
\newtheorem{definition}{\noindent Definition}
\newtheorem{corollary}{\noindent Corollary}
\newtheorem{conjecture}{\noindent Conjecture}

\newtheorem{remark}{\noindent Remark}

\newcommand{\proof}[1]{{\bf Proof.} #1 \rule{2mm}{2mm} \medskip}

\begin{document}

\begin{titlepage}

\begin{center}
{\LARGE \bf Statistical properties of braid groups \\ in locally free
approximation

}
\bigskip

{\large A.M. Vershik $^{1,*}$, S. Nechaev $^{2,3}$ and R. Bikbov
$^{3}$}
\medskip

{\sl
$^{1}$ St. Petersburg Branch of Steklov Mathematical Institute, Fontanka
27, \\ 119011 St.Petersburg, Russia
\medskip

$^{2}$ UMR 8626, CNRS-Universit\'e Paris XI, LPTMS, Bat.100,
Universit\'e Paris Sud, \\ 91405 Orsay Cedex, France
\medskip

$^{3}$ L.D. Landau Institute for Theoretical Physics, Kosygin str. 2, \\
117940 Moscow, Russia
}
\bigskip

\end{center}


\begin{abstract}

Statistical and probabilistic characteristics of locally free group with
growing number of generators are defined and their application to statistics
of braid groups is given.
\end{abstract}

\tableofcontents
\bigskip





\hrule
\bigskip

\noindent
$^{*}$ Also at: Institut des Hautes \'Etudes Scientifiques, 35, route de 
Chartres, \\ \hspace*{18mm} 91440 Bures-sur-Yvette, France

\end{titlepage}

\section{Introduction}

The theory of knots, since works of Artin, Alexander, Markov (jr.) and
others, is traditionally connected to the braid group $B_n$ (see for review
\cite{9}). Jones polynomials \cite{11} and a whole subsequent progress
in knot theory is
based on representations of braid groups and Hecke algebras. But alongside
with a known problem of construction of topological invariants of knots and
links, it should be
noted a number of similar, but much less investigated problems.
We mean calculation of a probability of a knot formation in a given
homotopic type with respect to a uniform measure on a set of all closed
nonselfintersecting contours of a fixed length on a lattice. The
given problem, known as a "knot entropy" problem so far has not an adequate
mathematical apparatus and since recently has been studied mainly by
numerical methods \cite{1}. Nevertheless it is clear that the mentioned
set of problems is connected by large to random walks on noncommutative
groups.

The last years have been marked by occurrence of a great number of problems
of a physical origin, dealing with probabilistic processes on
noncommutative groups. Let us mention some examples, important for us. First
of all, that are the problems of statistics and topology of chain molecules
and related statistical problems of knots (see, for example, \cite{2}), as
well as the classical problems of the random matrix theory and localization
phenomena \cite{3}.

In a set of works \cite{4} a problem concerning calculation of a
mathematical expectation of a "complexity" of randomly generated knot had
been formulated, where the degree of
any known algebraic invariant (polynomials of Jones, Alexander, HOMFLY and
other) had been served for the characteristics of complexity.

As to the theory of random walk on braid groups, a particular number of
results devoted to the investigation conditional limiting behavior of
Brownian bridges on the group $B_3$ is known only \cite{5}. Thus, neither
Poisson boundary, nor explicit expression of harmonic functions for braid
groups are so far found.

In the present paper we consider statistical properties of locally
free and braid groups following the idea of the first author (A.V.)
developed in \cite{20} and extended in the papers
\cite{13,6,6a}. For study of braid groups we introduce the concept of so
named {\it locally free groups}, which is the particular case of {\it
local groups} in sense of works \cite{13,12}. This concept gives us very
useful tool for bilateral approximations for the number of nonequivalent
words in the braid groups and
semi--groups. Very important and apparently rather new aspect of this
problem consists in passing to the limit $n\to \infty$ in the group
$B_n$; just this limit is considered in our work. We found rather
unexpected stabilization of various statistical characteristics of the
local groups in this limit.

In \cite{12} the systematical approach to computation of various numerical 
characteristics
of countable groups is proposed. The essence of this approach deals
with simultaneous consideration of three numerical constants, properly
characterizing the logarithmic volume, as well as the entropy and the
escape (the drift) of the uniform random walk on the group. It happens
that the inequality which relates these three constants proved in
\cite{13,12} reflects the deep statistical properties of local groups.
The evaluation of number of words,
entropy and other statistical characteristics for the locally free groups
permits one to estimate the appropriate characteristics for the braid
groups. So, the locally free groups play the role of the
approximant to the braid groups. At the same time the study of locally
free groups has appeared to be useful for other models of statistical
physics, connected with problems of directed growth, theory of parallel
computations and etc. (see \cite{8} for review).

\section{Main definitions and statement of a problem}

\begin{definition}[Artin braid group] The braid group $B_{n+1}$ of $n+1$
"strings" has $n$ of generators $\{\sigma_1,\ldots,\sigma_n\}$ with the
following relations:
\be \label{1}
\left\{\begin{array}{ll}
\sigma_i\sigma_{i+1}\sigma_i =
\sigma_{i+1}\sigma_i\sigma_{i+1} & \qquad (1\le i<n) \\
\sigma_i\sigma_j=\sigma_j\sigma_i & \qquad (|i-j|\ge 2) \\
\end{array} \right.
\ee
\end{definition}
In connection with general definitions of the braid group, there
exists an extensive literature---see \cite{9}; for the last work on the
normal forms of words, we shall note \cite{10}.

\begin{figure}[ht]
\centerline{\epsfig{file=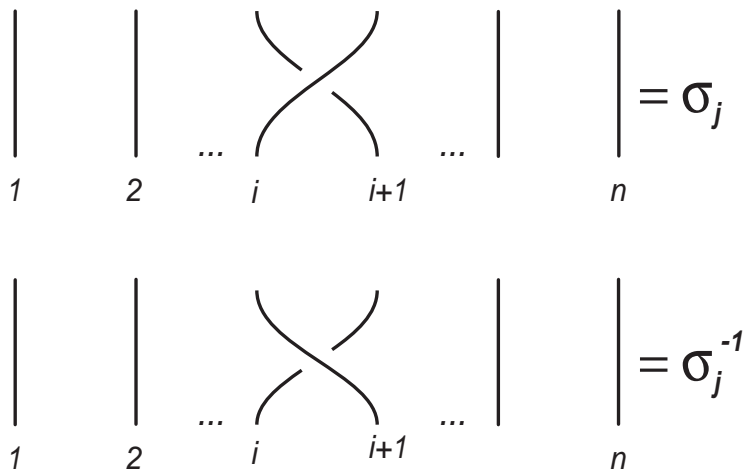,width=10cm}}
\caption{Graphic representation of generators of braid group $B_{n+1}$}
\label{fig1}
\end{figure}

The element of the braid group $B_n$ is set by a word in the alphabet
$\{\sigma_1,\ldots,\sigma_n; \break
\sigma_1^{-1},\ldots,\sigma_n^{-1}\}$---see fig.\ref{fig1}

By {\it length $N$ of a record} of a braid we call just a length of a word
in a given record of the braid, and by irreducible length (or simple {\it
length}) -- the minimal length of a word, in which the given braid can
be written. The irreducible length can be also viewed as a distance from
the unity on the Cayley graph of the group. Graphically the braid it is
represented by a set of strings, going from above downwards in accordance
with growth of a braid length.

Closed braid is received by gluing the "top" and "bottom" free ends on
a cylinder. Closed braid defines a link (in particular, a knot). Homotopic
type of the link can be described in terms of algebraic
characteristics of a braid \cite{11}.

Now we define a concept of a locally free group, which is a special case of
noncommutative local group \cite{13,12}.

\begin{definition}[Locally free group] By a locally free group
${\cal LF}_{n+1}$ with $n$ generators \break $\{f_1,\ldots,f_n\}$ we denote
a group, determined by following relations:
\be \label{2}
f_j f_k=f_k f_j \qquad \mbox{for all} \qquad |j-k|\ge 2, \quad
\{j,k\}=1,...n
\ee
\end{definition}
Each pair of neighboring generators $(f_j, f_{j\pm 1})$ produces a free
subgroup of a group ${\cal LF}_n$. Similarly the locally free semi--group
 ${\cal LF}_n^{+}$ is defined.

The concept, equivalent to the concept of locally free semi--group ${\cal
LF}_n^+$ has occurred earlier in work \cite{14}, devoted to the
investigation of combinatoric properties of substitutions of sequences and
so named "partially commutative monoids" (see \cite{8} and references
there). Especially productive becomes the geometrical interpretation of
monoids in a form of a "heap", offered by G.X. Viennot and connected with
various questions of statistics of directed growth and parallel
computations. The case of a group (instead of semi--group) introduces a
number of additional complifications to the model of a heap and apparently
has not been considered in the literature. We touch it in more details
below.

Obviously, the braid group $B_n$ is the factor--group of a locally free
free group ${\cal LF}_n$, since it is received by introducing the
Yang--Baxter (braid) relations. It has been realized also, that the group
${\cal LF}_n$ is simultaneously the subgroup of the braid group.

\begin{lemma} Consider a subgroup $\overline{B}_n$ of the group $B_n$,
over squares of generators $B_n:\; \overline{B}_n=\left<\overline
{\sigma}_1,..., \overline{\sigma}_{n-1}|
\overline{\sigma}_i=\sigma_i^2,\;
i=1,..., n\right>$. The correspondence $\overline{\sigma}_i \leftrightarrow
f_i$ sets the isomorphism of the groups $\overline{B}_n$ and ${\cal LF}_n$.
\end{lemma}

\proof{ For the proof of this fact it is sufficient to check that between
the generators $\overline{\sigma}_i=\sigma_i^2$ and $\overline{\sigma}_{i+1}
= \sigma_{i+1}^2$ there are no any nontrivial relations. Thus, it is
sufficient to restrict ourselves to consideration of a group $B_3$, or to be
more precise, of its subgroup $\overline{B}_3$. Consider the Burau
representation
$$
\sigma_1(t)=\left(\begin{array}{cc} 1 & 0 \\ t & -t \end{array}\right);
\qquad \sigma_2(t)=\left(\begin{array}{cc} -t & 1 \\ 0 & 1 \end{array}
\right),
$$
being the exact representation of $B_3$ over $C[t]$. It is obvious that
\be \label{sl}
\overline{\sigma}_1=\sigma_1^2(t)=\left(\begin{array}{cc} 1 & 0 \\
t-t^2 & t^2 \end{array}\right); \qquad \overline{\sigma}_2=\sigma_2^2(t)=
\left(\begin{array}{cc} t^2 & -t+1 \\ 0 & 1 \end{array} \right)
\ee
Putting $t=-1$, we see that (\ref{sl}) is reduced up to
$$
f_1=\sigma_1^2(-1)=\left(\begin{array}{cc} 1 & 0 \\ -2 & 1 \end{array}
\right); \qquad f_2=\sigma_2^2(-1)=\left(\begin{array}{cc} 1 & 2 \\ 0 & 1
\end{array}\right),
$$
which are the generators of free group, $\Gamma_2$.
}

It should be noted, that the matrices $\left(1\;2 \atop 0\;1\right)$ and
$\left(\;1\; \;0 \atop -2\; 1\right)$ are the generators of a free group.
This fact was proved apparently for the first time by I. Sanov \cite{16}.

\begin{figure}[ht]
\centerline{\epsfig{file=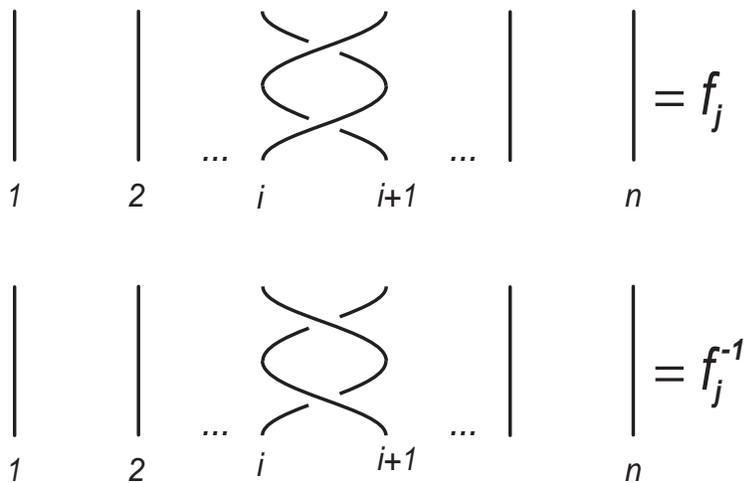,width=10cm}}
\caption{Graphic representation of generators of locally--free group ${\cal
LF}_{n+1}$}
\label{fig2}
\end{figure}

\begin{corollary}
Locally free group ${\cal LF}_n$ is simultaneously over-- and subgroup
of the braid group $B_n$.
\end{corollary}
This consequence will be hereafter used for transmitting the estimates from
the locally free group to the braid group. The geometrical interpretation
of the group ${\cal LF}_{n+1}$ is shown in fig.\ref{fig2}.

Let us introduce also a concept of a locally free group of restricted order.
\begin{definition}
We call the group ${\cal LF}_n^{(r)}$ with generators $\{f_1,\ldots, f_n\}$
and relations $(f_i)^r=1\; (i=1,...,n) \; f_i f_j = f_j f_i,\; |i-j| \ge 2$
a locally free group {\it of restricted order $r$}. The neighboring
generators $f_i, f_{i\pm 1}$ form the free product of ${\Z}/r{\Z}$.

In addition, another semi--group ${\cal LFP}_n^{+}$ emerges while imposing
"projective" relations $(f_i)^2=f_i$.
\end{definition}

Let us formulate the main problems concerning the determination of
asymptotic characteristics of locally free and similar groups.
\bigskip

\noindent {\bf 1. Asymptotics of number of words in a group.} Let $G$ be
the group with fixed framing $\{f_1,..., f_n\}$. The definition following
hereafter makes sense for any groups with fixed and finite set of
generators. Denote by $K(g)$ and call {\it the length} $K(g)$ the minimal
length of word $g$, written in terms of generators $\{g_1,...,g_n;
g_1^{-1},..., g_n^{-1}\}$. The length defines the metrics (the metrics
of words \cite{17}) on the group.

Denote by $V(G,K)$ the number of elements of group $G$ of length $K$.
\begin{definition}
Call $v(G,K)$ the {\it logarithmic volume of a group} with respect to the
given group $G$:
\be \label{ob}
v(G)=\lim_{K\to\infty} \frac {\log V(G,K)}{K}
\ee
\end{definition}
where the limit exists---see \cite{12}. We call $G$ as the group of
exponential growth if $v>0$.

In Section \ref{raz1} we investigate the asymptotic behavior of $v(B_n)$,
$v({\cal LF}_n)$ at $K\to \infty$, in the limit $n={\rm const}\gg 1$.
\bigskip

\noindent {\bf 2. Random walk and average drift on a group.} Consider
the (right--hand side) random walk on any group $G$ with fixed framing
$\{g_1,..., g_n; g_1^{-1},..., g_n^{-1}\}$, i.e. regard the Markov chain
with following transition probabilities: the word $w$ transforms into
$w\, g^{\pm 1}$ with the probability $\frac{1}{2n};\; i=1,..., n$. Similarly
one can build a left--hand Markov chain.

Let $L(G,N)$ be a mathematical expectation of a length of a random word,
obtained after $N$ steps of random walk on the group $G$.
\begin{definition}
Call $l(G,N)$ the {\it drift} on the group $G$ (see \cite{12}):
\be \label{snos}
l(G)=\limsup_{N\to\infty} \frac{L(G,N)}{N}
\ee
\end{definition}
Thus, the drift is the average speed of a flow to the infinity in the
metrics of words.

In Section \ref{raz2} we calculate the drift $l({\cal LF}_n)$ on the locally
free group and its limit for $n={\rm const}\gg 1$.
\medskip

\noindent {\bf 3. Entropy of a random walk on a group.} Let $\mu^N$ be the
$N$--time convolution of a uniform measure $\mu$ on generators
$\{f_1,..., f_n; f_1^{-1},..., f_n^{-1}\}$.
\begin{definition}
The {\it entropy} (see \cite{26,25,20,12}) of random walk on a group with
respect to $\mu$ refers to as
\be
\label{entr} h(G)=\lim_{N\to\infty} \frac{H\left(\mu^N\right)}{N} =
\inf_{N}\frac{H(\mu^N)}{N}
\ee
where $H(\nu)=-\sum\limits_{x\in {\rm supp}\; \nu}\nu(x)\log \nu(x)$.
\end{definition}

Section \ref{raz3} is devoted to the computation of $h({\cal LF}_n)$ in the
limit $n={\rm const} \gg 1$.

The question about simultaneous study of these three numerical
characteristics (volume, drift and entropy) is delivered by the first
author (A.V)---see \cite{12} and represents a serious and deep problem. In
particular, the desire to find the above defined characteristics for the
braid group motivates our consideration of locally free and similar groups.

The basic relation connecting these three numerical quantities is as follows
\be \label{3const}
v\,l\ge h
\ee
(see \cite{12} and also \cite{23}). For the free group this inequality is
reduced to the identity (\cite{12}).

\section{Asymptotics of number of words in locally free group}
\label{raz1}

Here we find the asymptotics in $n\gg 1$ of the logarithmic volume
of the group ${\cal LF}_n$ (see also \cite{5,24}). Later on, in Section
\ref{raz6} we use the results obtained here for the bilateral estimation
of the logarithmic volume of the braid group.

\begin{lemma}
Any element of length $K$ in the group ${\cal LF}_{n+1}$ can be uniquely
written in the normal form
\be \label{2:norm}
W=\left(f_{\alpha_1}\right)^{m_1}\left(f_{\alpha_2}\right)^{m_2}\ldots
\left(f_{\alpha_s}\right)^{m_s},
\ee
where $\sum_{i=1}^s|m_i|=K\;(m_i\neq 0\; \forall\; i;\; 1\le s\le K)$, and
indices $\alpha_1,...,\alpha_K$ satisfy the following conditions
\begin{itemize}
\item[(i)] If $\alpha_i=1$ then $\alpha_{i+1}=2,..., n$;
\item[(ii)] If $\alpha_i=k$ ($2\le k\le n-1$) then $\alpha_{i+1}=
k-1,k+1,...,n$;
\item[(iii)] If $\alpha_i=n$ then $\alpha_{i+1}=n-1$.
\end{itemize}
\end{lemma}

\proof{ The proof directly follows from the definition of commutation
relations in the group ${\cal LF}_{n+1}$.
}

It is easy to see, that the local rules (i)-(iii) define a Markov chain on
a set of the  generator numbers $\{1,...,n\}$ with $n\times
n$--dimensional transition matrix $\widehat{T}_n$
\be \label{matrix}
\widehat{T}_n=
\left(\begin{array}{ccccccc}
0 & 1 & 1 & 1 & \ldots & 1 & 1 \\
1 & 0 & 1 & 1 & \ldots & 1 & 1 \\
0 & 1 & 0 & 1 & \ldots & 1 & 1 \\
0 & 0 & 1 & 0 & \ldots & 1 & 1 \\
\vdots & \vdots & \vdots & \vdots & \ddots & \vdots & \vdots \\
0 & 0 & 0 & 0 & \ldots & 0 & 1 \\
0 & 0 & 0 & 0 & \ldots & 1 & 0
\end{array}\right)
\ee

\begin{theorem}
1) The number $V(n,K)$ of group elements ${\cal LF}_n$ of length $K$ is
equal to
\be \label{V}
V(n,K)=\sum_{s=1}^{K} 2^s\left(K-1 \atop s-1\right)
\th_n(s)=C\; \frac{2^n}{n^3} 7^K\Big(1+o_K(1)\Big)
\ee
where $\th_n(s)$ is the sum over all {\it various} sequences of $s$
generators ($1\le s\le K$), satisfying the rules (i)--(iii), and $C$
is the
numerical constant.

2) In the limit of infinite number of generators ($n\to \infty$) the
logarithmic volume of a locally free group is equal to
\be \label{vlog}
v=\lim_{K\to \infty} \frac{\log V(n,K)}{K}= \log 7
\ee
i.e. $v$ asymptotically corresponds to the logarithmic volume of a free
group with four generators.
\end{theorem}

\proof{ 1. The value $V(n,K)$ can be represented in the form
\be \label{2:trace}
V(n,K)=\sum_{s=1}^{K}\th_n(s)\mathop{{\sum}'}_{\{m_1,\ldots, m_s\}}
\delta\left[\sum_{i=1}^s|m_i|-K\right],
\ee
where the second sum gives the number of all representations of a word of
a given irreducible length $K$ for the {\it fixed} sequence of indices
$\alpha_1, ...\alpha_s$; prime means that the sum does not contain
the terms with $m_i=0$ ($1\le i\le K$); and $\delta(x)$ is the Kronecker
$\delta$--function: $\delta(x)=1$ for $x=0$ and $\delta(x)=0$ for $x\neq 0$.

The value $\th_n (s)$ makes sense of a partition function, determined as
follows:
\be \label{rs}
\th_n(s)=\left<{\bf v}\left[\widehat{T}_n\right]^{s-1}
{\bf v}\right>
\ee
where
\be \label{vec}
{\bf v}=(\;\overbrace{1\; 1\; \ldots\; 1}^{n}\;)
\ee
The remaining sum in expression (\ref{2:trace}) is independent on
$\th_n(s)$ and can be easily computed:
\be \label{2:perm}
\mathop{{\sum}'}_{\{m_1,\ldots,m_s\}} \delta\left[\sum_{i=1}^s
|m_i|-K\right] = 2^s\; C_{K-1}^{s-1}
\ee
where $\disp C_{K-1}^{s-1}=\frac{(K-1)!}{(s-1)!(K-s)!}$.

Substituting (\ref{2:perm}) and (\ref{rs}) in (\ref{2:trace}), we arrive at
the first statement of the theorem:
\be \label{fin}
V(n,K)=\sum_{s=1}^{K} 2^s\; C_{K-1}^{s-1} \; \th_n(s)
=  2{\bf v} (2\widehat{T}_n + \widehat{I})^{K-1}{\bf v}
\ee
where $\widehat{I}$ there is the identity matrix.

Our approach to the calculation of $\th_n(m)$ is based on a consideration
of a "correlation function" $\th_n(x,x_0,s)$, which determines the number
of various sequences of $s$ generators, satisfying the rules (i), (ii),
(iii), beginning with $f_{x_0}$ and finishing with $f_x$. Using the
representation (\ref{matrix}) it is not hard to write down an evolution
equation in "time" $s$ for the function $\th(x,s)\equiv \th_n(x,x_0,s)$:
\be \label{5}
\th(x,s+1)=\th(x-1,s)+\sum_{y=x+1}^{n}\th(y,s),
\ee
where
$$
\th_n(x,x_0,s)=\left<{\bf x}_0\left[\widehat{T}_n\right]^{s-1}
{\bf x}\right>,
$$
and
$$
{\bf x}_0=(\; \underbrace{\overbrace{0\; \ldots\; 0\;
1}^{x_0}\;0\;\ldots\;
0}_{n}); \quad {\bf x}=(\; \underbrace{\overbrace{0\; \ldots\; 0\;
1}^{x}\;0\;\ldots\; 0}_{n})
$$
The equation (\ref{5}) should be completed by initial and boundary conditions
on a segment $[0,n+1]$:
\be \label{6}
\left\{\begin{array}{l}
\th(x,1)=\delta_{x,x_0} \\
\th(0,s)=\th(n+1,s)=0
\end{array}\right.
\ee

Passing from (\ref{5}) to the local difference equation, we arrive at the
following boundary problem on a segment $[0,n+1]$:
\be \label{theta}
\left\{\begin{array}{l}
\disp \theta(x+1,s+1)-\theta(x,s+1)=
\theta(x,s)-\theta(x-1,s)-\theta(x+1,s) \\
\disp \theta(x,1)=\delta_{x_0,x} \\
\disp \theta(0,s)=\theta(n+1,s)=0.
\end{array}\right.
\ee
For a generating function defined as
\be
Z(x,s)=\sum_{s=1}^{\infty} p^{s-1} \th(x,s),
\ee
and using (\ref{theta}), we get an equation on $Z(x,p)$:
\be
\left\{\begin{array}{l}
(p+1)Z(x+1,p)-(p+1)Z(x,p) + pZ(x-1,p)= \delta_{x,x_0-1}-\delta_{x,x_0}
\\
Z(0,p)=Z(n+1,p)=0
\end{array}\right.
\ee
This equation can be symmetrised via substitution $Z(x,p)=A^x\varphi(x,p)$,
where $A=\sqrt{\frac{p}{p+1}}$, what results in a boundary problem:
\be \label{7}
\left\{\begin{array}{l}
\varphi(x+1,p)-\frac{1}{A}\varphi(x,p)+\varphi(x-1,p)= \disp
\frac{A^{-x}}{\sqrt{p(p+1)}} (\delta_{x,x_0-1}-\delta_{x,x_0}) \\
\varphi(0,p)=\varphi(n+1,p)=0.
\end{array}\right.
\ee

The Fourier $\sin$--transform $\disp f(k,p)=\sum_{x=1}^{n} \varphi (x,p)
\sin\frac{\pi kx}{n+1}$ allows one to rewrite (\ref{7}) as follows
$$
\left(2\cos\frac{\pi k}{n+1} -\frac{1}{A}\right)f(k,p)=
\frac{A^{-x_0}}{\sqrt{p(p+1)}}\left(A\sin\frac{\pi k(x_0-1)}
{n+1}-\sin\frac{\pi kx_0}{n+1}\right),
$$
and the final expression for $Z(x,p)$ reads
$$
Z(x,p)=\frac{2}{(n+1)(p+1)}
\left(\frac{p}{p+1}\right)^\frac{x-x_0}{2}\sum_{k=1}^{n}
\frac{ \sqrt{\frac{p+1}{p}}\sin\frac{\pi kx_0}{n+1} -\sin\frac{\pi
k(x_0-1)}{n+1}}{\sqrt{ \frac{p+1}{p}}-2\cos\frac{\pi k}{n+1}}\sin
\frac{\pi kx}{n+1}.
$$

It is convenient to express the function $\th(x,s)$ via the contour
integral
$$
\th(x,s)=\frac{1}{2\pi i}\oint\limits_{C_0}
\frac{Z(x,p)}{p^{s}}dp,
$$
where the contour $C_0$ surrounds an origin $p=0$ and lies in the
regularity area of the function $Z(p)\equiv Z(x,p)$. Hence,
$$
\th(x,s)=-\sum\limits_{p_k} \mbox{Res}
\left(\frac{Z(x,p_k)}{p_k^{s}}\right),
$$
where $p_k$ are the poles of the function $Z(p)$:
\be \label{poles}
p_k=\frac{1}{4\cos^2\frac{\pi k}{n+1}-1}
\ee

We are interested only in the asymptotic behavior of the function
$\th_n(x,x_0,s)$ at $s\gg 1$, which is determined by the poles nearest to
the origin, $p_1=p_{n}=\frac{1}{3}$ for $n \gg 1$. Thus, we have:
\be
\th_n(x,x_0,s)=\frac{4}{n+1}\sin\frac{\pi
(x_0+1)}{n+1} \sin\frac{\pi x}{n+1}\;2^{x_0-x}\;3^{s-1}
\ee

In order to find $\th_n(s)$ it is necessary to sum up $\th(x,x_0,s)$ over
all $x,x_0$: $\disp \th_n(s)=\sum_{x,x_0=1}^{n} \th_n(x,x_0,s)$. That gives
at $n={\rm const} \gg 1$ the following expression:
\be \label{18}
\th_n(s)=\frac{16\pi^2}{\log^4 (2/e)}\;\frac{2^n}{n^3}\;3^{s-1}\,
\Big(1+o(1)\Big)
\ee
Finally, using (\ref{V}), we have the following asymptotic expression (with
$n\gg 1$) for the total number of nonequivalent irreducible words of length
$K\gg 1$
\be \label{vv1}
V(n,K)=\frac{32\pi^2}{\log^4 (2/e)}\;\frac{2^n}{n^3}\;7^{K-1}\,
\Big(1+o_K(1)\Big)
\ee
Using the definition (\ref{ob}), we obtain the statement 2) of the theorem:
$$
\frac{\log V(n,K)}{K}\sim \log 7 + \frac{n-3\log n+C}{K} = \log 7 + o_K(1)
$$
where $C=\frac{32\pi^2}{\log^4 (2/e)}$. So,
$$
\lim_{K\to\infty} \frac{\log V(n,K)}{K}=\log 7
$$

Thus, for large number of generators, the logarithmic volume of the group
${\cal LF}_n$ saturates tending to the value $\log 7$, which corresponds to
the free group with four generators.
}

\begin{remark}
The expression (\ref{rs}) can be rewritten in the following compact form
$$
V(n,K)=\left<{\bf v}\left[\widehat{M}_n\right]^K {\bf v}\right>
$$
where the vector ${\bf v}$ is defined by Eq.(\ref{vec}) and the matrix
$\widehat{M}_n$ is as follows:
$$
\widehat{M}_n=
\left(\begin{array}{ccc}
\begin{array}{ccccc}
0 & 1 & 1 & 1 & \ldots \\
1 & 0 & 1 & 1 & \ldots \\
0 & 1 & 0 & 1 & \ldots \\
0 & 0 & 1 & 0 & \ldots \\
\vdots & \vdots & \vdots & \vdots & \ddots
\end{array} & \vline &
\begin{array}{ccccc}
0 & 1 & 1 & 1 & \ldots \\
1 & 0 & 1 & 1 & \ldots \\
0 & 1 & 0 & 1 & \ldots \\
0 & 0 & 1 & 0 & \ldots \\
\vdots & \vdots & \vdots & \vdots & \ddots
\end{array} \\ \hline \vspace{-0.3in} \\
\begin{array}{ccccc} \vspace{-0.15in} \\
0 & 1 & 1 & 1 & \ldots \\
1 & 0 & 1 & 1 & \ldots \\
0 & 1 & 0 & 1 & \ldots \\
0 & 0 & 1 & 0 & \ldots \\
\vdots & \vdots & \vdots & \vdots & \ddots
\end{array} & \vline &
\begin{array}{ccccc} \vspace{-0.15in} \\
0 & 1 & 1 & 1 & \ldots \\
1 & 0 & 1 & 1 & \ldots \\
0 & 1 & 0 & 1 & \ldots \\
0 & 0 & 1 & 0 & \ldots \\
\vdots & \vdots & \vdots & \vdots & \ddots
\end{array}
\end{array}\right)
$$
\end{remark}

\begin{remark}
Let us note, that the poles $p_k$ (see (\ref{poles})) have the
single--valued correspondence with the eigenvalues $\lambda$ of the matrix
(\ref{matrix}): $\lambda_k=p_k^{-1}$. In turn, the determinant of the matrix
$a_n(\lambda)=\det\left(\widehat{T}_n-\lambda\widehat{I}\right)$ satisfies
the recursion relation (see \cite{24})
\be \label{co:2}
\left\{\begin{array}{l}
a_n(\lambda)=-(\lambda+1)a_{n-1}(\lambda)-(\lambda+1)a_{n-2}(\lambda)
\\
a_0(\lambda)=1 \\
a_1(\lambda)=-\lambda
\end{array}\right.
\ee
At $\lambda>-1$ the solution of (\ref{co:2}) is expressed via Chebyshev
polynomials of the second kind:
\be \label{co:8}
a_n(\lambda) = (-1)^n(\lambda+1)^{\frac{n-1}{2}}\,
{\cal U}_{n+1}(\cos\vartheta) = (-1)^n(\lambda+1)^{\frac{n-1}{2}}\,
\frac{\sin(n+2)\vartheta}{\sin\vartheta}
\ee
where
\be \label{co:7}
\cos\vartheta=\frac{\sqrt{\lambda+1}}{2} \qquad
\left(0<\vartheta<\frac{\pi}{2}\right)
\ee
\end{remark}

\begin{corollary}
1. The number $V(n,K|{\cal LF}_n^+)$ of elements of length $K$ of a
semi--group ${\cal LF}_n^+$ is equal
$$
V(n,K|{\cal LF}_n^+)=\sum_{s=1}^{K}\left(K-1 \atop s-1\right)\th_n(s)=
C\; \frac{2^n}{n^3} 4^K(1+o_K(1))
$$
In the limit $n\to \infty$ the logarithmic volume of the locally free
semi--group reads
$$
v({\cal LF}_n^+)=\log 4+o_K(1)
$$

2. The number $V(n,K|{\cal LF}_n^{(r)})$ of elements of length $K$ of the
locally free group ${\cal LF}_n^{(r)}$ of local degree $r$ has the following
expression
$$
V(n,K|{\cal LF}_n^{(r)})=\sum_{s=1}^{K}{\cal N}_r(K,s)\th_n(s)
$$
where
$$
\begin{array}{lll}
\disp {\cal N}_r(K,s) & = & \disp \mathop{{\sum}'}_{\scri
\shortstack{$\{m_1,\ldots,m_s\}$ \\ $|\mbox{mod $r$}|$}}
\delta\left[\sum_{i=1}^s |m_i|-K\right] \medskip \\
& = & \left\{\begin{array}{ll} \delta_{K,s}, & r=2 \medskip \\
\disp \frac{1}{2\pi i}\oint\limits_{(C_1)} dz\; z^{-1-K+s}
\left(\frac{2-z^{m-1}+z^m}{1-z}\right)^s, & r=2m\; (m\ge 2) \medskip \\
\disp \frac{1}{2\pi i}\oint\limits_{(C_1)} dz\; z^{-1-K+s}
\left(\frac{2-z^m}{1-z}\right)^s, & r=2m+1\; (m\ge 1)
\end{array}\right.
\end{array}
$$
and the contour $C_1$ surround an origin of the complex plane $z$.

In particular, for $r=2,3,4$ in the limit $n\to \infty$ the logarithmic
volume of the locally free group of the local degree $r$ is:
$$
\begin{array}{l}
v({\cal LF}_n^{(2)})\to \log 3 \medskip \\
v({\cal LF}_n^{(3)})\to \log 6 \medskip \\
v({\cal LF}_n^{(4)})\to \log\left(3+2\sqrt3\right) \approx \log 6.464
\medskip \\ \disp \lim_{s\to\infty} \lim_{n\to \infty}
v({\cal LF}_n^{(s)})= \log 7
\end{array}
$$
The similar computations allows one to obtain the logarithmic volume
$v({\cal LFP}_n^{+})$ of the locally free semigroup with projective
relation:
$$
v({\cal LFP}_n^{+})= \log 3
$$
\end{corollary}

\section{Random walk on locally free group: the drift}
\label{raz2}
The computation of the drift of the random walk on the locally free group
presented below generalize the appropriate results for the free group.

Remind that a symmetric random walk on a free group $\Gamma_n$ with $n$
generators can be viewed as a cross product of a nonsymmetric
random walk on a half--line ${\Z}^+$ and a layer over $N\in{\Z}^+$ giving
a set of all words of length $N$ with the uniform distribution. The
transiton probabilities in a base are:
$$
N\to \left\{\begin{array}{ll}
N+1 & \mbox{with the probability $\frac{2n-1}{2n}$} \medskip \\
N-1 & \mbox{with the probability $\frac{1}{2n}$}
\end{array}\right.
$$
Thus, the mathematical expectation of a word's length after $N$ steps reads
$$
E \sum_{i=1}^N \xi_i = \sum_{i=1}^N E \xi_i = N \, E \xi_i =
N\left(\frac{2n-1}{2n} - \frac{1}{2n} \right)=N\, \frac{n-1}{n}
$$
and hence the drift is
$$
l=\lim_{N\to\infty} \frac{1}{N} E\sum_{i=1}^N \xi_i = \frac{n-1}{n}
$$
For example, for the group $\Gamma_2$ with two generators ($n=2$) the drift
is equal to $\frac{1}{2}$.

To compute the drift of the random walk on the group ${\cal LF}_n$ one
should understand in more details the structure of the normal form of elements
of ${\cal LF}_n$.

\begin{figure}[ht]
\centerline{\epsfig{file=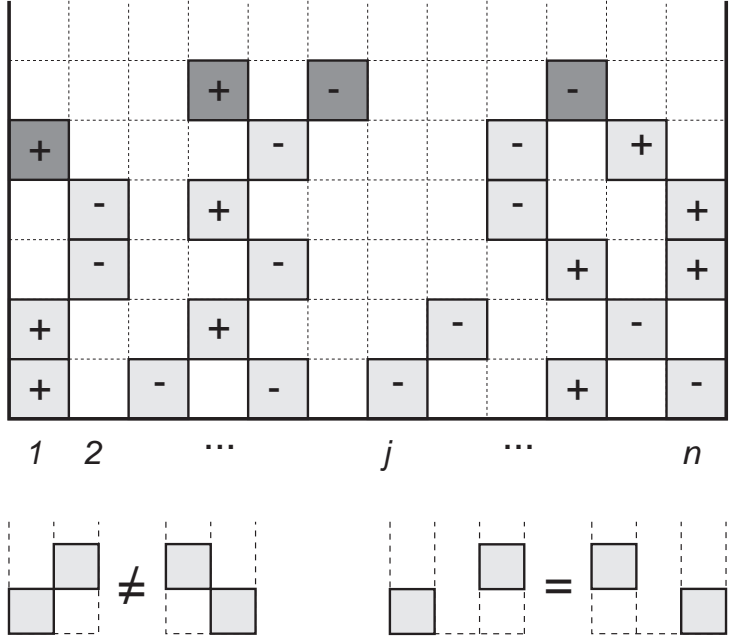,width=10cm}}
\caption{Typical configuration of a colored heap. Elements of a roof are
shown by filled squares.}
\label{fig3}
\end{figure}

It is helpful to use a geometrical interpretation of indicated concepts
following the ideas of G.X. Viennot (see \cite{8} for review), arised in 
connection
with the theory of partially commutative monoids \cite{14}. We imagine a
word in the group ${\cal LF}_n$ as a finite configuration of cells in a set
named herafter as a {\it heap} (a {\it colored heap})---see fig.\ref{fig3}.
\bigskip

Namely, we consider the strip $P={\bf n}\times{\Z}^{+}\subset{\Z}^2\;
\left({\bf n}=\{1,2,..., n\};{\Z}^+= \{0,1,...,\}\right)$ as a subset of
the lattice ${\Z}^2$.

\begin{definition}
I. We call as a {\it heap} the finite set $S\in P$, satisfying the
conditions:

1) In a horizontal line the cells from the set $S$ can not be neighbors;

2) Each cell from $S$, not standing in the first horizontal line has at
least one cell in the previous horizontal line, touching it. (The touching
of cells means that the horizontal coordinates of such cells differ not
more, than by 1.)

II. Let each cell from $S$ has two colors $(+,-)$. We shall assume, that
besides the conditions 1 and 2 the following one is fulfiled:

3) In one and the same vertical row $i=1,...,n$ the cells of different
color can not be the neighbors.

In the last case $S$ refers to as a colored heap.
\end{definition}

The set of heaps with number of cells $K\ge 0$ is denoted by $H_K^{(n)}$
(while $H^{(n)}=\mathop{\cup}\limits_{K=0}^{\infty}H_K^{(n)}$ ($H_0^{(n)}$
is an empty heap). The concept of a heap had been introduced and investigated
by Viennot \cite{8} in connection with combinatoric problems of
partially commutative monoids of Cartier and Foata \cite{14} and so-called
"directed lattice animals" considered in \cite{27,8}. Denote by $CH_K^{(n)}$
the set of colored heaps with the number of cells $K\ge 0$, thus
$CH^{(n)}=\mathop{\cup} \limits_{K=0}^{\infty} CH_K^{(n)}$. As far as we
know, the colored heaps have not been considered so far.

\begin{definition}
The numbered heap is a heap, whose cells are enumerated by natural
numbers so, that the the enumeration is monotone, i.e. if two cells touch
each others, the top cell has larger number.
\end{definition}

\begin{lemma}
There exists a bijection $\tau$  between the set of words of a locally free
group ${\cal LF}_n$ and the set of colored heaps $CH^{(n)}$, for which the
one-to-one correspondence between the elements (words) of ${\cal LF}_n$ and
colored heaps $CH^{(n)}$ is established.
\end{lemma}

For the semi--group ${\cal LF}_n^+$ the same statement is true with the
replacement $CH^{(n)}\to H^{(n)}$ and for this case the interpretation
was given in the work \cite{8}.

\proof{ By induction. The unity of the group corresponds to an empty heap.
Let $\tau$ be determined for words of length $\le K$. Compare to a word
of length $K+1$ a colored heap, which is received by adding the element
$g_i^{\pm 1}$ in an $i$'s vertical row to already existing colored heap,
i.e. put a cell so, that the conditions 1 and 2 of Lemma 3 are satisfied.
If directly under a new cell there was already a cell with the same
coordinate $i$ and the opposite color, these cells cancel. It is easy to
see, that if two words are equal, the corresponding (colored) heaps coincide.

\begin{figure}[ht]
\centerline{\epsfig{file=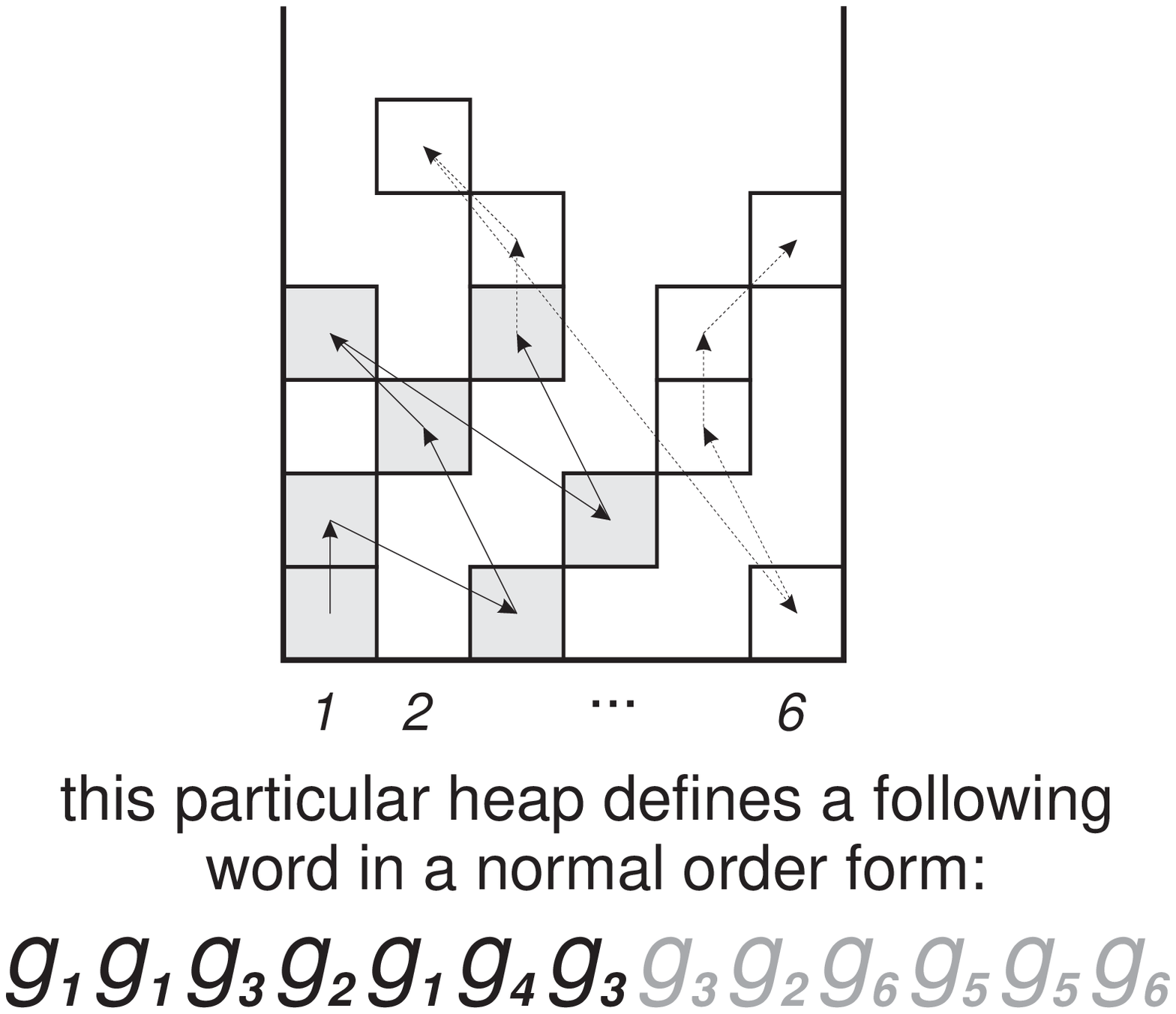,width=10cm}}
\caption{Example of construction of normally ordered word by given heap.
The configuration of some current cluster is shown by filled cells.}
\label{fig3a}
\end{figure}

Let us show now that any (colored) numbered heap is uniquely associated
with some word in ${\cal LF}_n$. Namely, we construct an algorithm which
sets a word in the normal order by some numbered (colored) heap:

1. Denote the most left--hand cell at bottom as the cell No.1 corresponding
to the first letter in the normal order form for a given heap. For
definiteness assume that this cell is disposed in a $j$'s column. The cell
No.2 is a cell located in a column $k$ ($k\le j$) as close as possible to
the cell No.1. Now we search for cells left--most close to cell No.2 and so
on... Continuing such enumeration we get a part of a heap called "cluster".

2. If there are no more cells satisfying the rule 1, we continue numbering
with the most bottom cell which is closest right--hand neighbor to the
given cluster {\it such that this new cell leaves the roof of the cluster
without changes}. This new cell is added to the cluster and enumeration is
continued recursively.

As a result get numbering corresponding the normal form of a given word.
Thus, Lemma 3 is proved.
}

\begin{remark}
There is an analogy between heaps and Young diagrams, as well as between
numbered heaps and Young tables \cite{12}.
\end{remark}

Dynamics of the words' growth, i.e. the random walk on ${\cal LF}_n$
(${\cal LF}_n^+$) acquires a following obvious geometrical sense: it
is a Markov chain with the states taken from the set of colored heaps
(or just heaps). The transitions consist in addition of cells (in a
view of conditions 1--3 of the Definition 7) with the probabilities
$\frac{1}{2n}$ for ${\cal LF}_n$ and $\frac{1}{n}$ for ${\cal LF}_n^+$.

For each element $w$ of the locally free group ${\cal LF}_n$, written
in the normal form we define the set of {\it achievable generators} $T(w)$
(see \cite{6a}):
$$
T(w)=\left\{i=1,...,n: K(w\,g_i)=K(w)-1\right\}
$$
where $K(w)$ is the word's length.

Achievable generator $g_i$ in the word $w$ refers to as such generator in
the representation $w=w_1\, g_i^m\, w_2$, that at multiplication of the
element $w$ from the right hand side by $g_i^{\pm 1}$, the new word $w'$
can be recorded as $w'=w_1\, g_i^{m\pm 1}\, w_2$, i.e. $w_2$ and $g_i$
commute. In particular, the generator $g_i^{m\pm 1}$ can be reduced, if
$m\pm 1=0$. It is easy to realize that achievable generators are such, that
can be reduced in one step of a random walk. Further we shall number
generators $g_i$ by the index $i$, and generators $g_i^{-1}$ by the index
$-i$.

The set $T(w)$ has following obvious properties:

(i) If $i\in T(w)$, then $-i\notin T(w)$ and if $-i\in T(w)$ then $i\notin
T(w)$;

(ii) If $i\in T(w)$ or $-i\in T(w)$ then $\pm (i-1), \pm (i+1)\notin
T(w)$.

The last property entails the inequality $\#T(w)\le \left[\frac{n+1}{2}
\right]$.

We continue the geometric interpretation of concepts and give the
visual description of the set $T(w)$ of achievable elements.

\begin{definition} \label{roof1}
We call as the {\it roof} of the heap (of the colored heap) the set $T(S)$
of those elements of a heap which have no upper neighbors in the same and
closest vertical rows.

In other words, some element belongs to the roof of the heap of $K$ elements
if after the removal of this element we get the heap of $K-1$ elements.
\end{definition}

\begin{remark}
In continuation of the similarity of heaps and Young diagrams, we can say
that the roof is analogous to the corners of the Young diagrams.
\end{remark}

\subsection{Mathematical expectation of the heap's roof}

In a geometrical interpretation described above the set $T(S)$ of
achievable elements is a set of such cells of a heap, removing of which for
one step leaves a configuration allowable (see the Definition \ref{roof1}).

The roof's basis\footnote{Hereafter, if is not stipulated especially, we
shall use the notation "roof" for a designation, both the roof as well as
the basis of a roof.} of a heap $S$ is the subset $T(S)=\{1,...,n\}$ of the
numbers of achievable generators. This subset, as it can be seen from the
properties (i)--(ii), satisfies the condition: if $(k_1, k_2) \in
T(S)$ then $|k_1-k_2|>1$. Denote by ${\cal T}_n$ the family of all such
subsets of the set $\{1,...,n\}$. In case of a colored heap the basis
consists of the painted in two colors $(+,-)$ subsets of $\{1,..., n\}$. We
denote these subsets by ${\cal T}_n^c$.

\begin{remark}
Let $w$ be the element of the group ${\cal LF}_n$ and $S$ be the
corresponding heap. Then $T_n(S)$ is exactly the set of achievable
generators.
\end{remark}

It is convenient to characterize ${\cal T}_n$ by a vector $(\ve_1,...,
\ve_n)$ with elements 0 and 1, where $\{\ve_r = 1\}\Leftrightarrow r\in T$.

\begin{lemma}
The power of the set ${\cal T}_n$ is equal to the Fibonacci number $F_n$ and
hence it grows as $\lam^n$, where $\lam=\frac{\sqrt{5}+1}{2}\approx 1.618$ is
the golden mean. The power of the set ${\cal T}_n^c$ is equal to $2^n$.
\end{lemma}

\proof{ The power of the set ${\cal T}_n$ is equal to the number of sequences
of elements $0$ and $1$ of length $n$, such that these sequences have not
the elements $1$ in succession, i.e. satisfy the recursion relation
\be
\begin{array}{l}
F_n=F_{n-1}+F_{n-2} \medskip \\ F_1=1; F_2=2,
\end{array}
\ee
which defines the Fibonacci sequence.

Similarly, the number of the elements of the set ${\cal T}_n^c$ satisfies
the recursion relation
\be \label{colfib}
\begin{array}{l}
F_n^c=F_{n-1}^c+2F_{n-2}^c \medskip \\
F_0^c=1; F_1^c=2; F_2^c=4
\end{array}
\ee

Actually, if the sequence ${\cal T}_n^c\subset F_n^c$ begins with 0,
the part remaining after removal of 0 is any sequence from $F_{n-1}^c$. If
${\cal T}_n^c$ begins with $1$, then by definition the 2nd element is 0.
Deleting these two elements (1 and following after it 0), we get a
sequence from $F_{n-2}^c$, Thus, the power of the set ${\cal T}_n^c$
satisfies recursion relation (\ref{colfib}), and consequently $F_n^c=2^n$.
}

Define the time--homogeneous Markov chain, the set of states of which in
any moment of a time are the sets $T\in{\cal T}$ and the transition
probabilities from the state $T$ to the state $T'$ are determined by the
time--independent rules. Let $T=\{\ve_1,...,\ve_n\};\; T'=\{\ve'_1,...,
\ve'_n\}$. Then the transition matrix is as follows. The transition
probability $T\to T'$ is nonzero and is equal to $\frac{1}{n}$
only for the cases when $\ve_i=\ve'_i$ for all $i$ except not more than
three consecutive numbers, say $(\ve_{r-1}, \ve_r, \ve_{r+1})$ and
$(\ve'_{r-1}, \ve'_r,\ve'_{r+1})$ and for these triples one of the
following conditions is satisfied:
\be \label{tripples}
\begin{array}{ll}
\mbox{If $\ve_{r-1}=\ve_{r+1}=1$} & \mbox{then $\ve'_r=1,
\ve'_{r-1}=\ve'_{r+1}=0$;} \medskip \\
\mbox{If $\ve_{r-1}=1,\ve_{r+1}=0$} & \mbox{then $\ve'_r=1; \ve'_{r-1}=
\ve'_{r+1}=0$;} \medskip \\
\mbox{If $\ve_{r+1}=1,\ve_{r-1}=0$} &\mbox{then $\ve'_r=1; \ve'_{r-1}=
\ve'_{r+1}=0$;} \medskip \\
\mbox{If $\ve_{r-1}=\ve_r=\ve_{r+1}=0$} &\mbox{then $\ve'_{r-1}=\ve'_{r+1}=0,
\ve'_r=1$.}
\end{array}
\ee
Thus, the Markov chain is determined on the set of states ${\cal T}_n$.

Later on we will be interested in the asymptotics of a mathematical
expectation of the size of a roof following the outline of the paper
\cite{6a}.

\begin{theorem}
The limit of the mathematical expectation of the number of achievable
generators for a random walk on the semi--group ${\cal LF}_n^{+}$ for
$n\gg 1$ (i.e. the limit of the mathematical expectation of the roof of
the heap) is
\be \label{13}
\lim_{N\to\infty} E\#T(w_N)=\frac{n}{3}
\ee
\end{theorem}

\proof{ Compute the mathematical expectation of a number of removable
(achievable) elements when we do not distinguish between
generators and their inverses, i.e. for the random walk on the semi--group
${\cal LF}_n^+$.

Represent the elements of the roof $T(w)$ (i.e. the number of achievable
generators) graphically by filled boxes on the diagram as it is shown
below:

\begin{center}
\unitlength=1.00mm
\special{em:linewidth 0.4pt}
\linethickness{0.4pt}
\begin{picture}(145.00,18.00)
\put(5.00,9.00){\line(0,1){2.00}}
\put(19.00,9.00){\line(0,1){2.00}}
\put(33.00,9.00){\line(0,1){2.00}}
\put(47.00,9.00){\line(0,1){2.00}}
\put(61.00,9.00){\line(0,1){2.00}}
\put(75.00,9.00){\line(0,1){2.00}}
\put(89.00,9.00){\line(0,1){2.00}}
\put(103.00,9.00){\line(0,1){2.00}}
\put(117.00,9.00){\line(0,1){2.00}}
\put(131.00,9.00){\line(0,1){2.00}}
\put(145.00,9.00){\line(0,1){2.00}}
\put(5.00,10.00){\line(1,0){140.00}}
\put(3.00,10.00){\rule{4.00\unitlength}{4.00\unitlength}}
\put(45.00,10.00){\rule{4.00\unitlength}{4.00\unitlength}}
\put(101.00,10.00){\rule{4.00\unitlength}{4.00\unitlength}}
\put(129.00,10.00){\rule{4.00\unitlength}{4.00\unitlength}}
\put(5.00,5.00){\makebox(0,0)[cc]{1}}
\put(47.00,5.00){\makebox(0,0)[cc]{4}}
\put(103.00,5.00){\makebox(0,0)[cc]{8}}
\put(131.00,5.00){\makebox(0,0)[cc]{10}}
\end{picture}
\end{center}
\vspace{-0.3in}

\noindent Here $n=11,\; \#T=4$.

Denote by $h_j=k_j-k_{j-1}-1$ the number of intervals of lengths $j$
between neighboring boxes or between a box and the edge of the
diagram.

Let $T$ consists of the set $\{k_1,...,k_s\}$. If the edge points $(1$ and
$n)$ do not belong to $T$, then $h_1=k_1,\; h_{s+1}=n-k_s-1$; if one or
both edge points belong to $T$, then $h_1=k_1-1,\; h_{s+1}=n-k_s$. For
example, if $k_1=1$ then $h_1=0$, or if $k_s=n$ then $h_s=0$. (On the above
diagram $h_1=0,\, h_2=1,\, h_3=3,\, h_4=1$). The numbers $h_j$ satisfy the
following relation, valid under neglecting the "boundary effects" at $\#T\gg
1,\; n\gg 1$:
\be \label{8}
\sum_j h_j=n-\#T
\ee
It is not hard to establish the rules according to which the diagram is
changed at such multiplication of $w$ by $g_r$ (or by $g_r^{-1}$), which
increases $\#T(w)$ by $1$: in $r$'s position appears a point while in
positions $(r-1)$ and/or $(r+1)$ the points (which were present) disappear.
Having in mind this rule, let us write the explicit expressions for the
1--step increment of a roof's length, $\Delta T(w)$, expressing it in terms
of $h_j(w)$ provided that the boundary points do not belong to $T$:
\be \label{loc}
\left\{\begin{array}{ll}
\Delta T(w)=+1 & \mbox{with the probability $\frac{1}{n}
\sum\limits_{j:h_j\ge 3}(h_j-2)$} \medskip \\
\Delta T(w)=0 & \mbox{with the probability $\frac{1}{n}\left(\#T+
2\#\{j:h_j\ge 2\}\right)$} \medskip \\
\Delta T(w)=-1 & \mbox{with the probability
$\frac{1}{n}\#\{j:h_j=1\}$}
\end{array}\right.
\ee
Summing (\ref{loc}), we obtain the conditional mathematical expectation of
the conditional probability of local reconstruction of a roof for the fixed
element $w$:
\be \label{loc2}
\begin{array}{lll}
E_w \Delta T & = & 1\cdot \frac{1}{n} \sum\limits_{j:h_j\ge 3} (h_j-2)+
0\cdot \left(\#T+ \frac{2}{n}\#\{j:h_j\ge 2\}\right)+ (-1)\cdot
\frac{1}{n}\#\{j:h_j=1\} \\
& = & \frac{1}{n} \sum\limits_j (h_j-2)= \frac{1}{n} \sum\limits_j
h_j - \frac{2}{n}\#T=1-\frac{3}{n}\#T(w)
\end{array}
\ee
Let us mention that depending on whether the boundary points belong or do
not belong to the set $T(w)$, the right--hand side of Eq.(\ref{loc2}) is
changed by terms which do not exceed than $\frac{4}{n}$. Therefore in the
large $n$ limit the expression (\ref{loc2}) is exact. In case of periodic
boundary conditions Eq.(\ref{loc2}) is exact for any finite values of $n$.

As far as our Markov chain has finite set of states and is ergodic, it has
the unique invariant measure. The Markov chain with this invariant measure
is stationary. So, the mathematical expectation $E[\Delta T(w)]$
over all elements $w$ with respect to the invariant measure exists and is
finite, therefore $E[\Delta T(w)]=0$. Thus, from the strong law of
large numbers (or, equivalently, from the individual ergodic theorem) it
follows that for the random walk on the semi--group we have Eq.(\ref{13})
for the mathematical expectation of the number of achievable elements (i.e.
the set of elements of the roof's elements).
}

Consider now the case of the group ${\cal LF}_n$. The distinction between
the semi--group ${\cal LF}_n^+$ (i.e. the heap) and the group ${\cal LF}_n$
(i.e. the colored heap) is due to the fact that for the random walk on the
group it is possible to reduce the word with the probability
$\frac{1}{2n}$. To account for that, we introduce the probabilities $p_w^+$
and $p_w^-$ to increase and to reduce the size of the roof $\#T(w)$ per unity
{\it under the condition of the word's length reduction}. The mathematical
expectation is a difference of conditional probabilities $p_w^{+}$ and
$p_w^{-}$ to change the value $\#T(w)$ per unity provided that reduction of a
word occurs. This difference should be added to the mathematical
expectation of the change of $\#T(w)$ in case of semi--group (\ref{tripples}):
\be \label{loc3}
E_w \Delta T = 1-\frac{3}{n}\#T(w) + \frac{p_w^{+} - p_w^{-}}{2n}\#T(w).
\ee
That gives (compare to (\ref{tripples})--(\ref{13}))
\be \label{13a}
E\#T=\frac{n}{3-\frac{1}{2}(p^{+}-p^{-})}=\frac{n}{3-\alpha}
\ee
where $p^{+}=Ep_w^{+};\; p^{-}=Ep_w^{-};$ and $\alpha=
\frac{1}{2}(p^{+}-p^{-})$.

On can easily realize that for some configurations of heaps
$w$ we could have $p^{+}-p^{-}\neq 0$ and in these cases the
mathematical expectation $E_w \Delta T$ for the group (colored heap) and
for the semi--group (heap) do not coincide. However, we believe, that at
$N\to \infty$ (i.e. in a stationary mode) $Ep_w^{+}=Ep_w^{-}$ and the
following hypothesis (expressed first by J. Desbois in \cite{6a}) is valid:
\begin{conjecture}
The mathematical expectation of a roof (a set of achievable elements) for
the heap (the locally free semi--group ${\cal LF}_n^+$) and for the colored
heap (the locally free group ${\cal LF}_n$) coincide at $n\gg 1$. Hence,
\be \label{loc4}
\lim_{N\to\infty} E\#T(w_N)=\frac{n}{3}
\ee
\end{conjecture}

The concept of a roof is the same for the heap (the semi--group) and for
the colored heap (the group), however the dynamics determined
in these two cases is distinct.

The random walk on the locally free semi--group (group) has been reduced to
a Markov dynamics of heaps (colored heaps). We have defined a new
dynamics---the dynamics of the roofs, Markovian in the case of
semi--group, by which the general dynamics is restored and which is
convenient for computations. In the case of the group this dynamics is not
Markovian anymore, but nevertheless enables to get some nontrivial
estimates.

\subsection{Drift as mathematical expectation of number of cells in the heap}

Let us compute now the change of a length of some fixed word $w$ for a
random walk on a group ${\cal LF}_n$. It is obvious, that for one step of
the random walk the length of a word can change by $\pm 1$. The
multiplication by a given generator, or by its inverse, occurs with the
probability $\frac{1}{2n}$ and thus, the conditional mathematical
expectation $E(w)$ to change a word's length is determined for a fixed
element $w$. Below we shall compute $E(w)$ and shall be convinced, that the
answer depends {\it only} on a size of a roof, i.e. on a size of a set
$\#T(w)$ of achievable generators $T(w)$.

Consider a fixed element $w$ of the group ${\cal LF}_n$ such that the set
of achievable generators $w$ is $\{1\le k_1<k_2<...<k_s<n\}$. Assume that
with the probability $\frac{1}{2n}$ the word $w$ is multiplied by a
generator $g_r$ or $g_r^{-1}$ (for definiteness let us choose $g_r$).
Denote the set of achievable generators of the element $w'=w\; g_r$ as
$T'\equiv T(w')$. Then the dynamics of the change of the set $T(w)$ is
settled by the following opportunities (compare to the above relations
(\ref{tripples})):

We have the following possibilities:

I. Provided that the the word is increased, i.e. $K'(w g_i)=K(w)+1$ the
dynamics of the roof is described by the relations (\ref{tripples}) valid for
the semigroup ${\cal LF}_n^{+}$;

II. Provided that the word is reduced, i.e. $K'(w g_i)=K(w)-1$, we have:
\be
\begin{array}{ll}
\mbox{If $\ve_r=1$} & \mbox{then ${\cal T}\to {\cal T}'\equiv {\cal T}^{-}$}
\end{array}
\ee
where ${\cal T}^{-}$ is the roof's configuration obtained by the
cancellation of one of the element of the roof ${\cal T}$ located in
position $r$. (This rule cannot be described in local terms).

The probability of a word's length reduction is $\frac{\#T}{2n}$, because
for each element of a roof there is a unique possibility to be reduced if
and only if at the following step the element inverse to the former one has
arrived. Accordingly, the probability to increase of a word's length is
$1-\frac{\#T}{2n}$, what follows from the possibility mentioned above to
change a word's length for one step by $\pm 1$. As a result, the mathematical
expectation of the total change of a word's length for one step of random
walk on the group
${\cal LF}_n$ is
\be \label{matozh}
E_{g_r}[K(w)-K(w\;g_r)]=
-\frac{E\#T(w)}{2n}+\left(1-\frac{E\#T(w)}{2n}\right) =
1-\frac{E\#T(w)}{n}
\ee
The indicated computation proves the following Lemma:
\begin{lemma} The conditional mathematical expectation of the word's length
$K(w)$ after $N$ steps of the random walk on the group ${\cal LF}_n$ for the
fixed last element $w$ is
$$
E_w K=N \left(1-\frac{E\#T(w)}{n}\right)
$$
hence the drift (i.e. the mathematical expectation of a normalized words'
length) is
$$
l=\lim_{N\to\infty}\frac{1}{N} E_w K = 1-\frac{E\#T(w)}{n}$$
\end{lemma}
Thus, for calculation of the drift it is sufficient to know the
mathematical expectation $E\#T(w)$ of the roof---see Eq.(\ref{loc3}).
\begin{theorem}
The mathematical expectation of the drift of a random walk on a locally
free group at $n\gg 1$ is
\be \label{loc4a}
l=1-\frac{E\#T(w)}{n}=1-\frac{1}{3-\alpha} = \frac{2-\alpha}{3-\alpha}
\ee
where $\alpha$ is defined in (\ref{13a}).
\end{theorem}

\begin{conjecture}
The mathematical expectation of the drift on the locally free group
at $n\gg 1$ is
\be \label{snos1}
l=\frac{2}{3}
\ee
\end{conjecture}
The Conjecture 2 is a direct consequence of the Conjecture 1 (J.
Desbois in \cite{6a}) but still it is not proved rigorously.

\section{Random walk on locally free group: entropy}
\label{raz3}

Let us remind, that the entropy $h$ of the random walk on the group $G$
according to the theorem similar to the Shannon's one and proved in
\cite{12,20,26,25} can be represented as follows (see the Definition
\ref{entr}):
$$
H(w)=-\lim_{N\to\infty}\frac{1}{N}\log\mu^N(g^N),
$$
where $w=g^N$; the limit in a right hand side is identical for almost
all trajectories $\{g^N\}_{N=1}^{\infty}$; $N$ is the number of the random
walk steps (i.e. the nonreduced word's length); $\mu^N$ is the $N$--time
convolution of a measure $\mu$.

In turn, the value of a measure $\mu^N$ on the words $w$ with $n$ generators
can be written in the following form
\be \label{mera}
\mu^N(w)=\frac{\#{\cal L}(w)}{n^N}
\ee
where $\#{\cal L}(w)$ is the number of various (dynamic) representations of the
element $w$ by words of length $N$ in the framing $\{g_1,..., g_n,
g_1^{-1},..., g_n^{-1}\}$. The value $\#{\cal L}(w)$ is the number of
different ways on the Cayley graph of the group, leading from the root
point of the graph (the unity of the group).

It is convenient to divide the problem of computation of the entropy of the
random walk on a group ${\cal LF}_n$ in two parts and to begin with the case
of the semi--group ${\cal LF}_n^+$, while for the group ${\cal LF}_n$ the
result will be the straightforward generalization of the corresponding
results for the semi--group ${\cal LF}_n^+$.

\bigskip

\subsection{Entropy of random walk on semi--group ${\cal LF}_n^+$}

As it has been found in the previous section during the study of the drift,
the dynamics of the increments of words (i.e. dynamics of the heap $S$) for
the random uniform addition of cells is uniquely determined by the dynamics
of the roof $T$ of the heap $H$. Moreover, we have found (see Eq.(\ref{13})),
that in the limit $N\to \infty$ and at $n\gg 1$ the mathematical expectation
of the roof's size $E\#T$, normalized by $n$ (i.e. the mathematical
expectation of the density of achievable elements) is $1/3$.

Let us prove the Lemma:
\begin{lemma}
The fluctuations of mathematical expectation of the roof for $n\gg 1$ have
the asymptotic behavior
$$
\frac{E\left|\#T^2-E(\#T)^2\right|}{E(\#T)^2}
\le\frac{{\rm const}}{n}
$$
where we have denoted
$$
E\#T\equiv \lim_{N\to\infty}E\#T(w_N)
$$
\end{lemma}

\proof{ Rewrite (\ref{13}) in the form
\be \label{13b}
E(\#T)^2\equiv \left[\lim_{N\to\infty}E\#T(w_N)\right]^2=\frac{n^2}{9}
\ee

Using Eqs.(\ref{loc})--(\ref{loc2}) for the probabilities of local
rearrangements of the roof we get the mathematical expectation of the
roof's fluctuations:
\be \label{entr1}
\begin{array}{lll}
E\Delta(\#T^2) & = & E\Delta\Big[(\#T')^2-(\#T)^2\Big] \medskip \\
& = & 1\cdot q_{+}\Big((\#T)^2-(\#T-1)^2\Big) + 0\cdot q_{0}+
(-1)\cdot
q_{-} \Big((\#T+1)^2-(\#T)^2\Big) \medskip \\
& = & 2(q_{+}-q_{-})\#T-(q_{+}+q_{-})
\end{array}
\ee
where
$$
q_{+}=\frac{1}{n}\sum_{j:h_j\ge 3}(h_j-2); \quad
q_{0}=\frac{1}{n}\left(\#T+ 2\#\{j:h_j\ge 2\}\right); \quad
q_{-}=\frac{1}{n}\#\{j:h_j=1\}
$$
Taking into account, that $q_{+}-q_{-}=1-\frac{3}{n} \#T$, we obtain from
(\ref{entr1}):
\be \label{entr2}
E\Delta(\#T^2)=E\left\{2\left(1-\frac{3}{n}\#T\right)\#T-
(q_{+}+q_{-})\right\}
\ee
For the invariant initial distribution it should be $E\Delta(\#T^2)=0$,
therefore the mathematical expectation of a square of the roof's size can
be received from the following relation
$$
E\Delta(\#T^2)=2E\#T-\frac{6}{n}E(\#T)^2-E(q_{+}+q_{-})=0
$$
whence we get
$$
E(\#T)^2=\frac{n}{3}E\#T-\frac{n}{6}E(q_{+}+q_{-})
$$
Estimating the mathematical expectation from above as $E(q_{+}+q_{-})<
{\rm const}$, we arrive at the equation:
$$
E(\#T)^2=\frac{n^2}{9}-\frac{{\rm const}\,n}{6}=\frac{n^2}{9}+o(n^2)
$$
Comparing the last expression with (\ref{13b}), we get the statement of
Lemma 6.
}

\begin{theorem}
The entropy of the random walk on the locally free semi--group ${\cal
LF}_n^{+}$ for $n\gg 1$ is
\be \label{entropy1}
h=\log 3 + O(n^{-1})
\ee
\end{theorem}

\proof{ Write the recursion relation of the change of the value $\#{\cal
L}(w)$ (see Eq.(\ref{mera}) at exception of one fixed cell of a
roof\footnote{Remind that the heap grows only by its roof.}, where
$w$ there is the element of the group, possessing at least one
representation by words of length $N$. Let $\#{\cal L}(w)$ and $\#{\cal
L}(w\setminus g_i)$ ($i\in T$) be the numbers of various representations of
the elements $w=g^N$ and $w'=g^N\setminus g_i$ correspondingly. Then is
fair the recursion relation \be \label{rek1} \#{\cal L}(g^N)=\sum_{g_i\in
T}\#{\cal L}(g^N\setminus g_i) \ee where the sum is taken over all elements
$g_i$ from the roof $T$.

Thus our problem is reduced to account for the number of representations
of a group element by generators in course of a random walk.

Using the definition (\ref{mera}) and Eq.(\ref{rek1}) we can write:
\be \label{rek2aa}
-\frac{1}{N}\log\mu^N(g^N) = -\frac{1}{N}\log\frac{\#{\cal L}(g^N)}{n^N}
\ee
Taking into account that the number of roofs of a given shape is equal to
the chronological ("time--ordered") multiplication of the roofs' lengths,
we arrive at the following equation:
\be \label{rek2ab}
-\frac{1}{N}\log\frac{\#{\cal L}(g^N)}{n^N}=
-\frac{1}{N}\log\frac{\disp \prod_{j=1}^{N}\#T_j}{n^N}=
-\frac{1}{N}\sum_{j=1}^{N}\log\frac{\#T_j}{n}
\ee
Defining the new variable
$$
\xi_j=\frac{3\#T_j}{n}-1
$$
we can rewrite (\ref{rek2ab}) in the form
\be \label{rek2a}
\begin{array}{c}
\disp \lim_{N\to\infty}-\frac{1}{N}\sum_{j=1}^N\log\frac{1+\xi_j}{3} =
\log 3 - \lim_{N\to\infty}\frac{1}{N}\sum_{j=1}^N\log(1+\xi_j) = \log 3 -
\lim_{N\to\infty} E\log(1+\xi_N) \\ \disp \le \log 3 - \lim_{N\to\infty}E\xi_N
+ \lim_{N\to\infty}\frac{1}{2}E\xi_N^2
\end{array}
\ee
where 
$$
\begin{array}{l}
\disp \lim_{N\to\infty}E\xi_N\equiv E\xi=\frac{3}{n}E\#T \\ 
\disp \lim_{N\to\infty}E\xi_N^2\equiv E\xi^2=\frac{9}{n^2}E\#T^2-
\frac{6}{n}E\#T+1,
\end{array}
$$ 
while $E\#T$ and $E\#T^2$ are the limits (for $N\to\infty$) of the 
mathematical expectations of the roof's size $\#T$ and of its square
 $\#T^2$ for fixed value of $n$.

The statement of Lemma 6 can be rewritten as
$$
\frac{E\left|\#T^2-\frac{n^2}{9}\right|}{\frac{n^2}{9}}=
E\left(\frac{3\#T}{n}\right)^2-1\le\frac{{\rm const}}{n}
$$
and in terms of the variable $\xi=\frac{3\#T}{n}-1$ the last inequality reads
$$
E(\xi+1)^2-1\le\frac{{\rm const}}{n}.
$$
Using the fact that $E\xi=0$ we arrive at the
conclusion that for $n\gg 1$
\be \label{rek2b}
E\xi^2\le \frac{{\rm const}}{n}
\ee

Substituting (\ref{rek2a}) for the estimate (\ref{rek2b}) we get
the desired statement of the Theorem:
\be \label{rek3}
h=\disp -\lim_{N\to\infty}\frac{1}{N}\,\log\mu^N(g^N) = \log 3+O(n^{-1})
\ee
what completes the proof.
}

\begin{theorem}
For the random walk on the locally free semi--group ${\cal LF}_n^{+}$ the
logarithmic volume $v$, the drift $l$ and the entropy $h$ satisfy at
$n\gg 1$ the strict inequality
$$
v\,l>h
$$
where $v\equiv v({\cal LF}_n^{+});\;l\equiv l({\cal LF}_n^{+});\;
h\equiv h({\cal LF}_n^{+})$.
\end{theorem}

\proof{

(i) From the Corrolary 2 of the Theorem 1 we have
$$
v({\cal LF}_n^{+})\to\log 4
$$

(ii) The drift of the random walk on ${\cal LF}_n^{+}$ is strictly equal to
1, i.e.
$$
l=1
$$

(iii) By Theorem 4 we have
$$
h({\cal LF}_n^{+})\to \log 3
$$

Comparing the values of $v$, $l$ and $h$, we get the strict inequality $v>h$
for the random walk on ${\cal LF}_n^{+}$.
}

\subsection{Entropy of random walk on group ${\cal LF}_n$}

\begin{theorem}
The entropy of the random walk on the locally free group ${\cal LF}_n$ at
$n\gg 1$ is
\be \label{entropy2}
h=\log (3-\alpha) + O(n^{-1})
\ee
where $\alpha=\frac{1}{2}(p^{+}-p^{-})$ (see \ref{13a}).
\end{theorem}

\proof{ In the case of the group the exact value of $E\#T$ is unknown as far
as $E\#T$ contains some unknown value $\alpha=\frac{1}{2}(E p_w^{+}-E
p_w^{-})$. Remind that $p_w^{+}$ and $p_w^{-}$ are the
probabilities of the change of the  value $\#T(w)$ by $\pm 1$ provided the
reduction of a word (see Eq.(\ref{13a})). Nevertheless, we can follow
directly the outline of the proof of the Theorem 4 with just a single
replacement $\xi_j\to \xi_j-\alpha$. As a result, in the limit $n\gg 1$ we
get (\ref{entropy2}).
}

For the group ${\cal LF}_n$ as well as in case of semi--group ${\cal
LF}_n^{+}$, the entropy $h$ and the drift $l$ of the random walk
are determined by the mathematical expectation of the roof's size $E\#T$.
Nevertheless in case of the group the numerical value of the mathematical
expectation of a colored heap's roof depends on the value $\alpha$. However
as far as our purpose is to prove that for locally free group in the limit of
infinite number of generators the strict inequality $l\,v<h$ holds, it is
sufficient to estimate appropriately the interval of change of $\alpha$.

\begin{theorem}
For the random walk on the locally free group ${\cal LF}_n^{+}$, the
logarithmic volume $v$, the drift $l$ and the entropy $h$ satisfy at
$n\gg 1$ the strict inequality
$$
v\;l> h
$$
where $v\equiv v({\cal LF}_n);\;l\equiv l({\cal LF}_n);\; h\equiv h({\cal
LF}_n)$.
\end{theorem}

\proof{ For the proof we shall use again the statements of the Theorems 1, 3,
and 4.

(i) By the Theorem 1:
\be \label{group1}
v({\cal LF}_n)\to\log 7
\ee

(ii) By the Theorem 3:
\be \label{group2}
l({\cal LF}_n)\to\frac{2-\alpha}{3-\alpha},
\ee

(iii) By the Theorem 4:
\be \label{group3}
h({\cal LF}_n)\to\log(3-\alpha).
\ee

By definition $\alpha= \frac{1}{2}(p^{+}-p^{-})$. Because $p^{+}+p^{-}=1$,
the following estimate is valid $|\alpha|<\frac{1}{2}$. Thus, the values of
the drift and the entropy lie within the interval
$$
\frac{3}{5}<l<\frac{5}{7}
$$
and
$$
\log \frac{5}{2} < h < \log \frac{7}{2}
$$

Define the discrepancy $\ve=l\,v-h$ and check that $\ve(\alpha)>0$ for all
values of $|\alpha|$ from the interval $|\alpha|<\frac{1}{2}$. Consider the
function
$$
\varepsilon(\alpha)=\frac{2-\alpha}{3-\alpha}\log 7 - \log(3-\alpha)
$$
Computing the derivative $\frac{d\ve(\alpha)}{d\alpha}$, one can easily
verify that on the interval $-\frac{1}{2}<\alpha<\frac{1}{2}$ the function
$\ve(\alpha)$ is strictly positive, hence $l\,v-h>0$. The Theorem is
proved.
}

\section{Conclusion}
\label{raz6}

Let us outline the applications of the results received above to the braid
groups and semi--groups.

\subsection{Some remarks about the relations to braid groups and
semi--groups}

As it has been mentioned already, the braid group $B_n$ is the
factor--group of the locally free group ${\cal LF}_n$ and simultaneously,
${\cal LF}_n$ is the subgroup of $B_n$. The same relations are valid for
the semi--group of positive braids $B_n^{+}$ and the locally free
semi--group ${\cal LF}_n^{+}$.

From here, as well as from the Consequence 1 of Lemma 1 and Theorem 1,
immediately follows the Theorem:
\begin{theorem}
The logarithmic volumes $v(B_n)$ and $v(B_n^{+})$ for $n\gg 1$ satisfy
the bilateral estimates
\be \label{otsenka}
\begin{array}{l}
\frac{1}{2}\log 7 < v(B_n) \le \log 7 \medskip \\
\log 2 < v(B_n^{+}) \le \log 4 \medskip
\end{array}
\ee
\end{theorem}

\proof{ The estimate from above is a direct consequence of the fact that
$B_n$ is the factor--group of the locally free group ${\cal LF}_n$. Thus,
$$
v_{B_n}\le v_{{\cal LF}_n}\equiv\log 7
$$

In order to obtain the estimate from below let us notice that the
embedding $\rho_n$ of ${\cal LF}_n$ in $B_n$ and of ${\cal LF}_n^{+}$ in
$B_n^{+}$ is realized via the identity $f_i\to \sigma_i^2$. Thus, in the
case of the group we have:
$$
V(\rho_n:\,{\cal LF}_n, K)\subset V(B_n,2K)
$$
and
$$
\frac{\log V(\rho_n\;{\cal LF}_n, K)}{K}\le \frac{V(B_n,2K)}{K},
$$
hence
$$
v_{{\cal LF}_n}\equiv\log 7 \le 2v_{B_n}
$$
therefore
$$
\frac{1}{2}v_{{\cal LF}_n}\le v_{B_n},
$$

Along the same line the case of the semi--group ${\cal LF}_n^{+}$ can be
treated.
}

Apparently, the upper estimate in Eq.(\ref{otsenka}) is closer to the true
value, than the lower one.

\begin{theorem}
The drift $l(B_n)$ on the braid group $B_n$ at $n\gg 1$ satisfies the
inequality
\be \label{otsenka2}
\frac{2-\alpha}{2(3-\alpha)} <l(B_n) \le \frac{2-\alpha}{3-\alpha}
\ee
\end{theorem}

\proof{ The bilateral estimate (\ref{otsenka2}) is also a direct
consequence of the fact that the the braid group is the factor--group of
the locally free group and in turn the locally free group is the subgroup
of the braid group. The value $\alpha$ has been defined above and varies in
the interval $-\frac{1}{2}<\alpha<\frac{1}{2}$.
}
\medskip

For the entropy of the random walk on the braid group the corresponding
bilateral estimates have not yet received.

\subsection{Physical interpretation of results}

\begin{figure}[ht]
\centerline{\epsfig{file=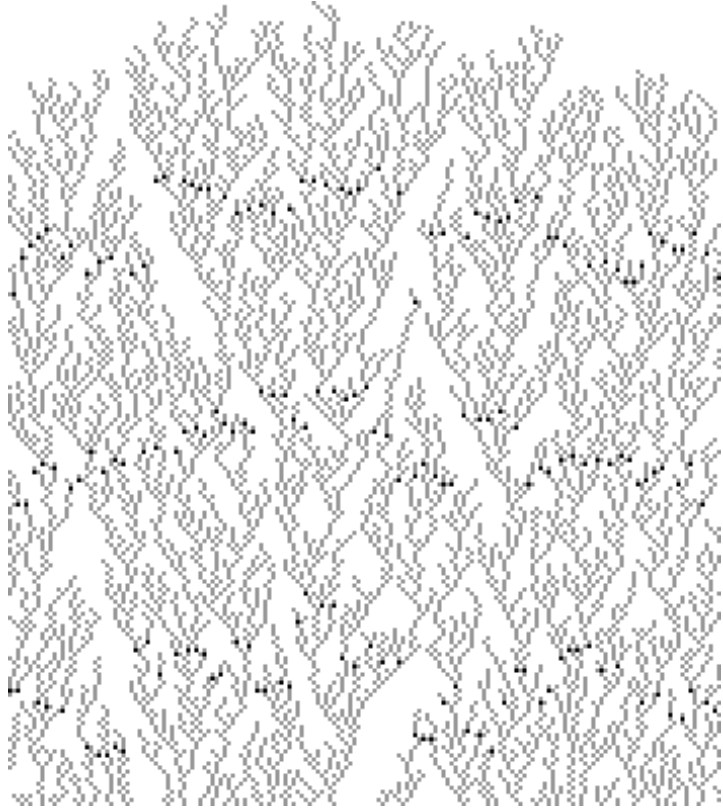,width=10cm}}
\caption{Typical configuration obtained in numerical simulations of the
uniform heap's growth.}
\label{fig4}
\end{figure}

Let us discuss briefly the physical sense of a strict inequality $lv>h$ for
the locally free group and for the ballistic deposition process---see
fig.\ref{fig4}.

The relation (\ref{mera}) permits one to estimate the number of various
dynamic representations of almost all (typical) elements $w$ by words of
length $N$ with respect to the uniform measure $\mu$:
$$
\#{\cal L}(w)\approx\exp(Nh)
$$
that for locally free group gives with the exponential accuracy
$$
\#{\cal L}(w)\approx(3-\alpha)^N
$$
Thus as stated already, the value $\#{\cal L}(w)$ is the {\it weighted}
(with the measure $\mu$) number of various states of the Cayley graph of the
locally free group, visited by a trajectory of random walk of length $N$.

On the other hand, the expression
$$
\#{\cal V}(w)\approx V^L=\exp(Nl\,v)
$$
gives the exponential estimate of the number of all different words, met
for a time of random walk on the group. For the locally free group the
value $\#{\cal V}(w)$ can be represented as follows
$$
\#{\cal V}(w)\approx7^{\frac{2-\alpha}{3-\alpha}N}
$$
where $|\alpha|<\frac{1}{2}$. In other words, $\#{\cal V}(L)$ is the {\it
complete} number of all different states of the Cayley graph of the locally
free group, located at a distance of typical drift $L$ of trajectory from
the root point of the graph.

The inequality
\be \label{ner1}
\#{\cal V}(w)\gg\#{\cal L}(w)
\ee
means, that the number of typical (on a measure) trajectories of the random
walk on the locally free group is exponentially small fraction of all
trajectories of the same length.

The inequality, similar to (\ref{ner1}) in case of the locally free
semi--group reads
\be \label{ner2}
\#{\cal V}^{+}(w)\gg\#{\cal L}^{+}(w)
\ee
where $\#{\cal V}^{+}(w)\approx 4^N$ is the volume of the locally free
semi--group ${\cal LF}_n^{+}$ for $n\gg 1$ and $\#{\cal L}^{+}(w)
\approx 3^N$ is the entropy of the random walk on ${\cal LF}_n^{+}$ in the
same limit.

The dynamically induced probabilistic measure on the group (semi--group), i.e.
the representation of words by the random walks on a group (semi--group),
essentially differs from the uniform (on the words) measure. This difference
is manifested in the exponential divergence of the two quantities
$\#{\cal V}(w)$ and $\#{\cal L}(w)$ ---see Eq.(\ref{ner1}) (the same
is valid for the semi--group described by Eq.(\ref{ner2})).

The inequality (\ref{ner2}) seems to be the origin of the fact that in the
numerical simulations of a random heap's growth (fig.\ref{fig4}) it is
observed a strong
divergence  between normalized mathematical expectation (averaged density) 
of the roof $\overline{\rho}_{\rm roof}=E\frac{\#T}{n}$ 
(where $\overline{\rho}_{\rm roof}=\frac{1}{3}$) and the
mathematical expectation (averaged density) of a whole heap 
$\overline{\rho}_{\rm heap}=\frac{N}{nH}$ 
(where $H$ is the maximal height of a heap).

The value of $H$, obtained in various computer experiments is evaluated as
$H\approx 4.05 \frac{N}{n}$ (for references see \cite{zhang}), what 
corresponds to the density $\overline{\rho}_{\rm heap}\approx
0.247$. The same value of the density is observed in average for any flat
horizontal section of a heap. A fact of essential numerical distinction
between $\overline{\rho}_{\rm roof}$ and $\overline{\rho}_{\rm heap}$ means 
that the roof of a 
heap has nontrivial fractal structure lying in a strip of nonzero's width. 
In fig.\ref{fig4} we have shown by black points few current configurations of
the roofs in course of the heap's growth. As it can be seen, the roof's
configurations are far from the flat ones and exhibit apparently the
nontrivial fractal behavior, which would be interesting to compare with the
continuous models of the surface growth described by Kardar--Parisi--Zhang
(KPZ) theory (see, for review \cite{zhang}).

\subsection*{Acknowledgments}

We would thank S. Fomin for pointing us to the connection between heaps and 
partially commutative monoids, as well as to G.X. Viennot, B. Derrida and 
A. Comtet for fruitful comments; S.N. highly appreciates deep suggestions 
made by J. Desbois (see in \cite{6a}). All authors are grateful to the
RFBR grant 99--01--17931 for partial support.


\end{document}